\documentclass[10pt]{article}

\usepackage{charter} 
\usepackage{fullpage}
\usepackage{titlesec}
\usepackage{amsmath,graphicx}
\usepackage{amssymb}
\usepackage{cite}
\usepackage{amsthm}
\usepackage{float}
\usepackage{color,subcaption}
\usepackage{bm}
\usepackage{bbm}
\usepackage{algorithm}
\usepackage{algpseudocode}
\usepackage{enumitem}

\titleformat*{\section}{\large\bfseries}
\titleformat*{\subsection}{\normalsize\bfseries}
\titleformat*{\subsubsection}{\small\bfseries}

\def\defn{\,\triangleq\,}
\def\sgn{\mathsf{sgn}}
\def\indfcn{\mathbbm{1}}

\def\kbm{{\bm{k}}}

\def\xbm{{\bm{x}}}
\def\ybm{{\bm{y}}}
\def\zbm{{\bm{z}}}
\def\Lbm{{\bm{L}}}
\def\Hbm{{\bm{H}}}
\def\thetabm{{\bm{\theta}}}
\def\sigmabm{{\bm{\sigma}}}

\def\Psf{\mathsf{P}}
\def\Tsf{\mathsf{T}}
\def\Gsf{\mathsf{G}}
\def\Psfhat{\widehat{\mathsf{P}}}
\def\Gsfhat{\widehat{\mathsf{G}}}

\def\xbmhat{{\widehat{\bm{x}}}}

\def\rhat{{\widehat{r}}}

\def\xbmdetil{\widetilde{\bm{x}}}

\def\xbmast{{\bm{x}^\ast}}

\def\R{\mathbb{R}}

\def\E{\mathbb{E}}
\def\P{\mathbb{P}}

\def\argmin{\mathop{\mathsf{arg\,min}}}
\def\signSGD{\mathsf{signSGD}}
\def\signProx{\mathsf{signProx}}
\def\SPGM{\mathsf{SPGM}}
\def\prox{\mathsf{prox}}

\theoremstyle{definition}

\newtheorem{definition}{Definition}
\newtheorem{theorem}{Theorem}

\newtheorem{assumption}{Assumption}

\begin{document}

\title{\Large \textsf{signProx}: One-Bit Proximal Algorithm for\\ Nonconvex Stochastic Optimization}

{\normalsize\author{Xiaojian Xu\(^1\)
and Ulugbek~S.~Kamilov\(^{1,2,}\thanks{This material is based upon work supported by the National Science Foundation under Grant No.~1813910.}\) \\
\emph{\small \(^1 \) Department of Computer Science and Engineering,~Washington University in St.~Louis, MO 63130, USA.}\\
\emph{\small \(^2 \) Department of Electrical and Systems Engineering,~Washington University in St.~Louis, MO 63130, USA.}\\
\small\emph{email}: \texttt{\{xiaojianxu,kamilov\}@wustl.edu}
}}

\markboth{\textsf{signProx}: One-Bit Proximal Algorithm for\\ Nonconvex Stochastic Optimization}%
{Kamilov: \textsf{signProx}: One-Bit Proximal Algorithm for\\ Nonconvex Stochastic Optimization}

\date{}
\maketitle 

\begin{abstract}
Stochastic gradient descent (SGD) is one of the most widely used optimization methods for parallel and distributed processing of large datasets. One of the key limitations of distributed SGD is the need to regularly communicate the gradients between different computation nodes. To reduce this communication bottleneck, recent work has considered a one-bit variant of SGD, where only the sign of each gradient element is used in optimization. In this paper, we extend this idea by proposing a stochastic variant of the proximal-gradient method that also uses one-bit per update element. We prove the theoretical convergence of the method for non-convex optimization under a set of explicit assumptions. Our results indicate that the compressed method can match the convergence rate of the uncompressed one, making the proposed method potentially appealing for distributed processing of large datasets.\end{abstract}


\section{Introduction}
\label{Sec:Intro}

Efficient processing of large datasets is a fundamental problem in modern signal processing. In many applications, the task can be formulated as an optimization problem of the form
\begin{equation}
\label{Eq:Minimization}
\xbmhat = \argmin_{\xbm \in \R^n} \left\{f(\xbm)\right\} \quad\text{with}\quad f(\xbm) = d(\xbm) + r(\xbm),
\end{equation}
where the data-fidelity term $d$ penalizes mismatch to the data and the regularizer $r$ enforces desirable properties in $\xbm$ such as sparsity or positivity. For a differentiable function $f$, the solution of~\eqref{Eq:Minimization} can be approximated iteratively with the classical gradient method~\cite{Boyd.Vandenberghe2004, Nesterov2004}
\begin{equation}
\xbm^t \leftarrow \xbm^{t-1} - \gamma \nabla f(\xbm^{t-1}),
\end{equation}
where $\gamma > 0$ is the step size. However, when $f$ consists of a large number of component functions
\begin{equation}
\label{Eq:AverageComponents}
f(\xbm) = \frac{1}{K} \sum_{k = 1}^K f_k(\xbm),
\end{equation}
the cost of computing the full gradient $\nabla f$ can become prohibitively expensive. In such cases, it is common to rely on the \emph{stochastic gradient descent (SGD)}~\cite{Robbins.Monro1951} that approximates the gradient at every iteration either with that of a single component $f_k$ or with an average of $B$ component gradients as
\begin{equation}
\label{Eq:StochGrad}
\xbm^t \leftarrow \xbm^{t-1} - \gamma \nabla \widehat{f}(\xbm^{t-1}) \quad\text{with}\quad \nabla \widehat{f}(\xbm) = \frac{1}{B}\sum_{b = 1}^B \nabla f_{k_b}(\xbm),
\end{equation}
where $k_1, \dots, k_B$ are independent random variables that are distributed uniformly over $[1 \dots K]$. 

A powerful feature of SGD is that it can be easily parallelized by splitting the computation of $B$ gradients over multiple compute nodes~\cite{Li.etal2014}. However, distributed SGD suffers from a significant communication overhead due the frequent gradient updates transmitted between the nodes. As the size of the gradient scales proportionally to the number of optimization parameters, it can reach hundreds of millions of variables for certain large-scale applications such as 3D imaging~\cite{Kamilov.etal2016} or deep learning~\cite{Huang.etal2017}. Motivated by this problem, recent work has considered a compressed SGD, where the algorithm compresses $\nabla \widehat{f}$ during optimization~\cite{Seide.etal2014, Alistarh.etal2017, Bernstein.etal2018}. A particularly simple variant of compressed SGD is $\mathsf{signSGD}$~\cite{Bernstein.etal2018}, which only keeps the sign of the stochastic gradient at every iteration
\begin{equation}
\xbm^t \leftarrow \xbm^{t-1} - \gamma \, \sgn\hspace{-0.2em}\left(\nabla \widehat{f}(\xbm^{t-1})\right),
\end{equation}
and hence compresses each stochastic gradient to a single bit. Remarkably, it was shown that under some conditions this simple scheme can match the convergence rate of uncompressed SGD~\cite{Bernstein.etal2018}.

While the current formulation of $\signSGD$ is both conceptually elegant and widely applicable, it does not take advantage of the recent progress in the proximal optimization theory~\cite{Parikh.Boyd2014}. In many applications, the regularizer $r$ in~\eqref{Eq:Minimization} consists of functions $r_k$ with easily computable proximals~\cite{Moreau1965}
\begin{equation}
\prox_{\gamma r_k}(\ybm) \defn \argmin_{\xbm \in \R^n}\left\{\frac{1}{2}\|\xbm-\ybm\|_2^2 + \gamma r_k(\xbm)\right\},
\end{equation}
which enables efficient optimization with a class of methods known as proximal algorithms. For example, two widely popular methods for large-scale optimization, FISTA~\cite{Figueiredo.Nowak2003, Bect.etal2004, Daubechies.etal2004, Bioucas-Dias.Figueiredo2007, Beck.Teboulle2009} and ADMM~\cite{Eckstein.Bertsekas1992, Ng.etal2010, Boyd.etal2011}, are both examples of proximal algorithms. 

In this paper, we propose a novel framework for one-bit stochastic optimization based on the proximal-gradient extension of $\signSGD$. Our method, called $\signProx$, is similar to $\signSGD$ in the sense that it also uses only one-bit per element of the update. In fact, we will see that under some conditions $\signProx$ is exactly equivalent to $\signSGD$. On the other hand, $\signProx$ also enables gradient-free optimization, and hence generalizes $\signSGD$ to problems with easily computable proximals. One of the key contribution of this paper is the theoretical analysis of $\signProx$ for nonconvex optimization under a set of transparent assumptions. Our analysis and simulations reveal that $\signProx$ can converge as fast or faster than the noncompressed algorithm, which makes compressed proximal optimization appealing for processing large datasets.


\begin{figure*}
\begin{minipage}[t]{.5\textwidth}
\begin{algorithm}[H]
\caption{$\textsf{SPGM}$}\label{alg:spgm}
\begin{algorithmic}[1]
\State \textbf{input: } $\xbm^0 \in \R^n$, $\gamma > 0$, and $B \geq 1$
\For{$t = 1, 2, \dots$}
\State sample a vector $\kbm$ with i.i.d. elements $k_b \sim \thetabm$
\State $\xbm^t \leftarrow \Psfhat_\kbm(\xbm^{t-1}) $
\EndFor\label{euclidendwhile}
\end{algorithmic}
\end{algorithm}%
\end{minipage}
\hspace{0.25em}
\begin{minipage}[t]{.5\textwidth}
\begin{algorithm}[H]
\caption{$\textsf{signProx}$}\label{alg:signProx}
\begin{algorithmic}[1]
\State \textbf{input: } $\xbm^0 \in \R^n$, $\gamma > 0$, and $B \geq 1$
\For{$t = 1, 2, \dots$}
\State sample a vector $\kbm$ with i.i.d. elements $k_b \sim \thetabm$
\State $\xbm^t \leftarrow \xbm^{t-1} - \gamma \, \sgn(\xbm^{t-1}-\Psfhat_\kbm(\xbm^{t-1}))$
\EndFor\label{euclidendwhile}
\end{algorithmic}
\end{algorithm}%
\end{minipage}
\end{figure*}

\section{Main results}
\label{Sec:Results}

In this section, we present our main results. We first introduce $\signProx$ and then follow up by analyzing its convergence.

\subsection{Compressed optimization using proximals}

The central building block of our algorithm is the following proximal-gradient mapping
\begin{equation}
\label{Eq:BasicPG}
\Psf_k(\xbm) \defn \prox_{\gamma r_k}\left(\xbm-\gamma \nabla d(\xbm)\right), \quad k \in [1, \dots, K]
\end{equation}
which computes a gradient-step with respect to the function $d$ and then evaluates the proximal with respect to another function $r_k$, both with a step-size $\gamma > 0$. Throughout this paper, we will assume that $d$ is a smooth, but possibly nonconvex, function. On the other hand, to have a well-defined proximal, we assume that $r_k$ are all closed, proper, and convex functions. We also define the following convex combination of mappings in~\eqref{Eq:BasicPG}
\begin{equation}
\label{Eq:FullPG1}
\Psf(\xbm) \defn \E[\Psf_k(\xbm)] = \sum_{k = 1}^K \theta_k \Psf_k(\xbm),
\end{equation}
where we can express the sum as an expectation since $\theta_k \geq 0$ and $\sum_{k = 1}^K \theta_k = 1$. It is known that a convex combination of proximals gives another proximal~\cite{Bauschke.etal2008, Yu2013}, which means that there exists a closed, proper, and convex function $r$ such that
\begin{equation}
\label{Eq:FullPG2}
\Psf(\xbm) = \prox_{\gamma r}(\xbm - \gamma \nabla d(\xbm)).
\end{equation}

Algorithm~\ref{alg:spgm} summarizes a stochastic alternative to the traditional proximal-gradient method~\cite{Beck.Teboulle2009b, Bertsekas2011}, which we will call $\SPGM$ in this paper. Instead of evaluating the full proximal-gradient step~\eqref{Eq:FullPG1}, which can be costly for large $K$, it computes an average proximal-gradient over a mini-batch of size $B$
\begin{equation}
\label{Eq:StochasticPG}
\Psfhat_\kbm(\xbm) \defn \frac{1}{B} \sum_{b = 1}^B \Psf_{k_b}(\xbm) = \prox_{\gamma \rhat}(\xbm - \gamma \nabla d(\xbm)),
\end{equation}
where each $k_b \in [1 \dots K]$ is sampled independently according to the probability distribution $\thetabm$ from~\eqref{Eq:FullPG1}. When the dependence of $\Psfhat_\kbm$ on $\kbm$ is clear from context, we will sometimes omit the subscript $\kbm$ from the notation as in $\Psfhat$. The second equality in~\eqref{Eq:StochasticPG} is due to the fact that $\Psfhat$ is a convex combination of proximals, and hence itself a valid proximal of some convex function $\rhat$.

Algorithm~\ref{alg:signProx} summarizes the main contribution of this paper: one-bit compressed version of $\SPGM$. Similarly to $\signSGD$~\cite{Bernstein.etal2018}, it requires only a single bit for updating each element of the iterate. However, the update direction at iteration $t$ is given by the sign of the quantity
$(\xbm^{t-1} - \Psfhat(\xbm^{t-1}))$. The choice of this direction is deliberate as it coincides with the gradient-mapping defined as follows.
\begin{definition}
For an objective $f(\xbm) = d(\xbm) + r(\xbm)$ and a step-size $\gamma > 0$, the \emph{gradient mapping} is defined as the operator
$$\Gsf(\xbm) \defn \frac{1}{\gamma}(\xbm - \Psf(\xbm)) = \frac{1}{\gamma}(\xbm - \prox_{\gamma r}(\xbm-\gamma\nabla d(\xbm))), \,\,\forall \xbm \in \R^n.$$
\end{definition}
\noindent
It is common to analyze the convergence of proximal algorithms using the gradient mapping, since $\Gsf(\xbmast) = 0$ \emph{if and only if} $\xbmast$ is the critical point of $f$~\cite{Beck2017}. Hence, $\signProx$ simply uses a one-bit approximation for elements of the gradient mapping at every iteration. 

We conclude this section by noting that $\signProx$ can also be seen as a generalization of $\signSGD$. Let $d = 0$ and define $r_k$ to be a linear approximation of $f_k$ around $\xbm^{t-1}$
\begin{equation}
\label{Eq:LinearApprox}
r_k(\xbm) = f_k(\xbm^{t-1}) + \nabla f_k(\xbm^{t-1})^\Tsf(\xbm-\xbm^{t-1}).
\end{equation}
Then, one can verify that the stochastic proximal-gradient iteration~\eqref{Eq:StochasticPG} reduces to the SGD iteration~\eqref{Eq:StochGrad}
\begin{equation}
\Psfhat(\xbm^{t-1}) = \frac{1}{B}\sum_{b = 1}^B \Psf_{k_b}(\xbm^{t-1})\ = \xbm^{t-1} - \gamma \nabla \widehat{f}(\xbm^{t-1}).
\end{equation}
Which means that the update of $\signProx$ will reduce to
$$
\xbm^t = \xbm^{t-1} - \gamma \sgn(\xbm^{t-1}-\Psfhat(\xbm^{t-1})) = \xbm^{t-1}-\gamma \sgn(\nabla \widehat{f}(\xbm^{t-1})).
$$
Hence, $\signSGD$ can be interpreted as Algorithm~\ref{alg:signProx} applied to a linear approximation of a function.

\subsection{Theoretical analysis}

We now discuss the convergence of $\SPGM$ and $\signProx$. Convergence of the stochastic proximal-gradient algorithms in the convex setting was analyzed by Bertsekas~\cite{Bertsekas2011}. Here, we focus on the case where $d$ is nonconvex. Our result for $\signProx$ extends the analysis of $\signSGD$ in~\cite{Bernstein.etal2018} using the theory of proximal optimization.
\begin{assumption} 
\label{As:Assumption1}
We analyze $\SPGM$ under the following assumptions:
\begin{enumerate}[label=(\alph*)]
\item The objective function $f$ has a finite minimum $f^\ast = f(\xbmast)$ attained at some $\xbmast \in \R^n$.
\item The function $d$ is differentiable and has a Lipschitz continuous gradient with a constant $L > 0$.
\item All functions $r_k$ are closed, proper, and convex. We also assume that they have Lipschitz continuous gradients with the same constant $L > 0$.
\item The proximal-gradient mappings have a bounded variance
$$\E\left[\|\Psf_k(\xbm)-\Psf(\xbm)\|^2\right] \leq \gamma^2 \sigma^2, \quad \forall \xbm \in \R^n,$$
for some constant $\sigma > 0$, where $\gamma > 0$ is the step-size.
\end{enumerate}
\end{assumption}
\noindent
All these are standard assumptions used in the analysis of stochastic optimization algorithms. The dependence of the variance (d) on $\gamma$ might seem surprising; however, note that this comes from the dependence of the proximal-gradient on $\gamma$, with $\gamma = 0$ implying that $\Psf_k(\xbm) = \Psf(\xbm)$ for all $k \in [1 \dots K]$ and $\xbm \in \R^n$.

\begin{theorem}
\label{Thm:SPGMconv}
Run $\SPGM$ for $T$ iterations under Assumption~\ref{As:Assumption1} with the step $\gamma = 1/(L\sqrt{T})$ and the mini-batch size $B = 1$. Then, we have that
$$\E\left[\frac{1}{T}\sum_{t = 1}^T \|\Gsf(\xbm^{t-1})\|_2^2\right] \leq \frac{1}{\sqrt{T}}\left[2L(f(\xbm^0)-f^\ast) + 3\sigma^2\right].$$
\end{theorem}

\noindent
The proof is given in Section~\ref{Seq:Proof1}. This establishes that SPGM converges to the critical point of the objective $f$.

Our analysis of $\signProx$ will need the following more elaborate set of assumptions.
\begin{assumption} 
\label{As:Assumption2}
We analyze $\signProx$ under the following assumptions:
\begin{enumerate}[label=(\alph*)]
\item The objective function $f$ has a finite minimum $f^\ast = f(\xbmast)$ attained at some $\xbmast \in \R^n$.
\item The function $d$ is differentiable and there exists a nonnegative vector $\Lbm \defn (L_1, \dots, L_n)$ such that
$$|\nabla d(\xbm)_i - \nabla d(\ybm)_i| \leq L_i |x_i - y_i|, \,\,\, \forall i \in [1\dots n], \forall \xbm, \ybm \in \R^n$$
\item All functions $r_k$ are closed, proper, and convex. We additionally assume that they all satisfy
$$|\nabla r_k(\xbm)_i - \nabla r_k(\ybm)_i| \leq L_i |x_i - y_i|, \,\,\, \forall i \in [1\dots n], \forall \xbm, \ybm \in \R^n.$$
\item The proximal-gradient mappings have a bounded variance
$$\E\left[(\Psf_k(\xbm)_i-\Psf(\xbm)_i)^2\right] \leq \gamma^2 \sigma_i^2, \,\,\, \forall i \in [1\dots n],  \forall \xbm \in \R^n,$$
for a positive $\sigmabm \defn (\sigma_1, \dots, \sigma_n)$, where $\gamma > 0$ is the step-size.
\end{enumerate}
\end{assumption}

\noindent
Note (b) and (d) lead to the standard assumption of Lipschitz continuity by defining a Lipschitz constant $L \defn \|\Lbm\|_\infty$. Similarly, the standard variance bound is recovered by setting $\sigma^2 = \|\sigmabm\|_2^2$. Also note that when the mini-batch size is $B > 1$, the variance bound is effectively reduced by $B$ for the mini-batch.

\begin{theorem}
\label{Thm:signProxconv}
Run $\signProx$ for $T$ iterations under Assumption~\ref{As:Assumption2} with the step $\gamma = 1/(2\|\Lbm\|_1\sqrt{T})$ and the mini-batch size $B = T$. Then, we have that
$$\E\left[\frac{1}{T}\sum_{t = 1}^T \|\Gsf(\xbm^{t-1})\|_1\right] \leq \frac{4}{\sqrt{T}}\left[\|\Lbm\|_1(f(\xbm^0)-f^\ast) + \|\sigmabm\|_1 + 1\right].$$
\end{theorem}

\noindent
The proof is given in Section~\ref{Seq:Proof2}. One can see that $\signProx$ has the same $\ell_1$-geometry as $\signSGD$, where the convergence rate depends on the $\ell_1$-norm of the gradient mapping, the stochasticity via $\sigmabm$, and the curvature via $\Lbm$. Surprisingly, it is possible for $\signProx$ to outperform $\SPGM$, when the gradient-mapping is dense but has a sparse set of extremely noisy components (see the detailed discussion for $\signSGD$ in~\cite{Bernstein.etal2018}). Our simulations in the next section will highlight this situation by comparing the relative performances of $\SPGM$ and $\signProx$ for nonconvex phase retrieval. Finally, to conclude this section, note that the theoretical analysis here was done for nonconvex functions $f$. It would be very interesting to see how convexity of $d$ can strengthen these convergence results. Note that the majority vote in $\signSGD$~\cite{Bernstein.etal2018} also directly applies to our algorithm; however, we omit its analysis due to the page limitation.


\section{Numerical illustration}
\label{Sec:Experiments}

\begin{figure}[t]
\begin{center}
\includegraphics[width=10cm]{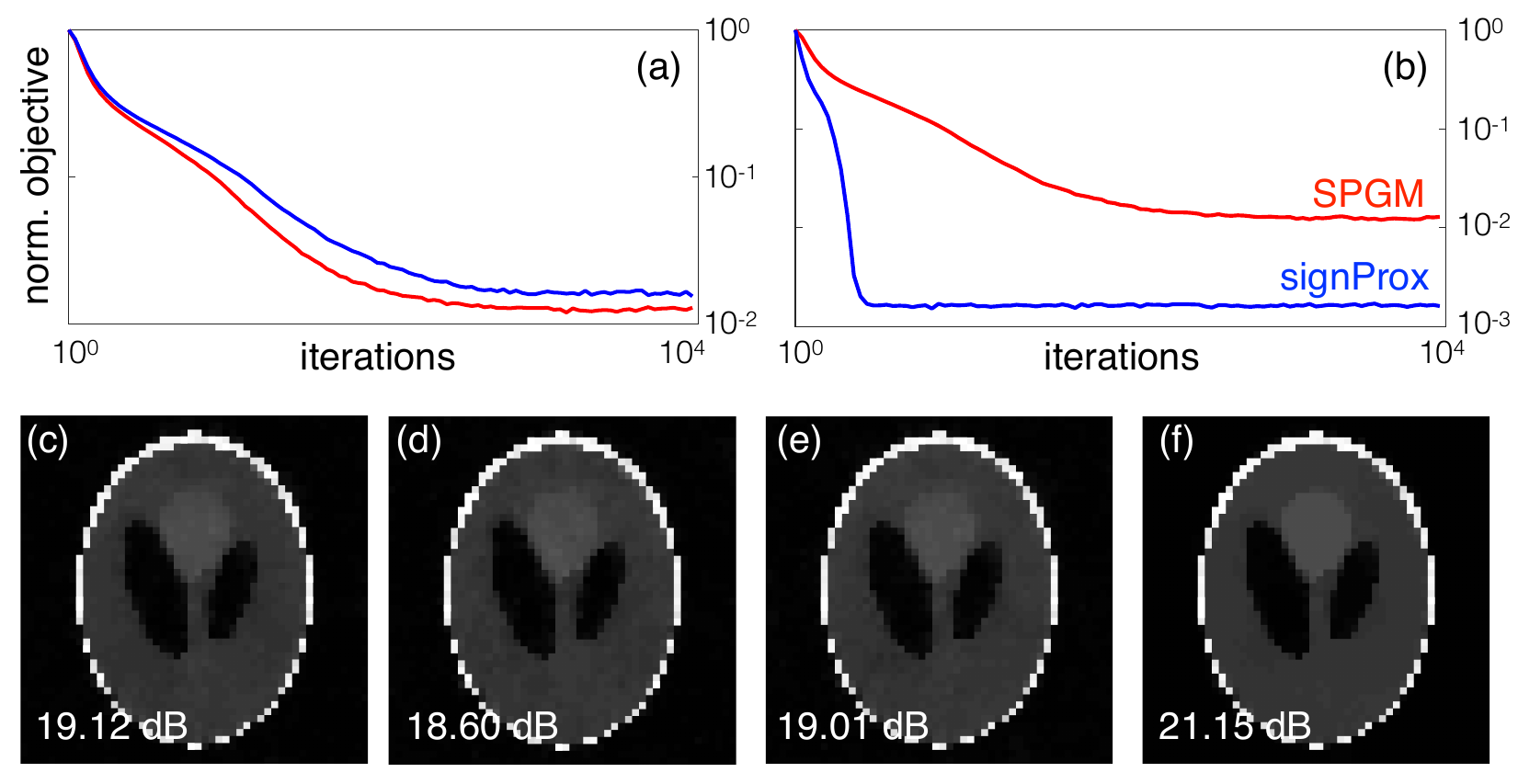}
\end{center}
\caption{Comparison of $\SPGM$ and $\signProx$ on the problem of generalized phase retrieval using TV regularization in two scenarios: with dense (a, c, d)  and sparse (b, e, f) stochasticity $\sigmabm$ of the proximal. The images (c, e) are the outputs of $\SPGM$, while (d, f) are the outputs of $\signProx$. This figure illustrates that in some settings $\signProx$ can converge similarly or even faster than $\SPGM$.}
\label{Fig:Scenario}
\end{figure}

We illustrate the relative performance of $\SPGM$ and $\signProx$ with a simple example. Consider the problem of generalized phase retrieval that was extensively considered in the literature~\cite{Candes.etal2013, Schniter.Rangan2012, Shechtman.etal2015}. When the signal is real, the goal is to reconstruct $\xbm \in \R^n$ given a set of nonlinear measurements $\ybm = |\zbm|^2$ with $\zbm = \Hbm\xbm \in \R^m$. This problem can be formulated as nonconvex optimization with a data-fidelity term $d(\xbm) = \frac{1}{2}\|\ybm-|\Hbm\xbm|^2\|_2^2$ and a sparsity-preserving regularizer such as total variation (TV)~\cite{Rudin.etal1992, Beck.Teboulle2009a, Kamilov2017}. Figure~\ref{Fig:Scenario} considers the reconstruction of a $50 \times 50$ \emph{Shepp-Logan phantom} from $m = 3000$ intensity measurements $\ybm$, where the measurement matrix $\Hbm$ is random with i.i.d.~$\mathcal{N}(0, 1/m)$ elements. We obtain a stochastic algorithm by using an unbiased estimate $\Psfhat$ for $\Psf$, obtained by adding a random noise to $\Psf$ distributed i.i.d~as $p(e) = \rho \mathcal{G}(e; \sigma_e^2) + (1-\rho)\delta(e)$, where $\mathcal{G}(e; \sigma_e^2)$ denotes the Gaussian pdf of variance $\sigma_e^2$ and $\delta(e)$ is the Dirac delta function. By setting $\rho \in (0, 1]$, we control the sparsity of the noise in the proximal-gradient, which directly corresponds to shaping the stochasticity $\sigmabm$ in Assumption~\ref{As:Assumption2}(d). Intuitively, when there is a sparse set of very noisy updates, the performance of $\SPGM$ will be dominated disproportionally by the noise, while the effect on $\signProx$ will be reduced as it discards the update amplitude. This is visible in Figure~\ref{Fig:Scenario}, where (a, c, d) correspond to the dense stochasticity scenario ($\rho = 1$) and (b, e, f) correspond to the sparse stochasticity scenario ($\rho = 0.1$). In both cases, the standard deviation of the noise is kept constant to the product of $\gamma$ and $\sigma = 0.1$. The step-size $\gamma$ is selected for the best performance. The convergence is quantified with the normalized objective ${(f(\xbm^t)-f^\ast)/(f(\xbm^0)-f^\ast)}$, where $f^\ast$ was obtained using the full TV reconstruction. One can observe that the convergence rate of $\SPGM$ is the same in both settings, while the convergence of $\signProx$ is faster when the stochasticity is sparse.


\section{Conclusion}
\label{Sec:Conclusion}

We have proposed a new $\signProx$ algorithm for stochastic optimization. The updates of $\signProx$ are compressed as each stochastic update contains a single-bit per element. We have proved the convergence of the method on nonconvex objectives under explicit assumptions. The future work will investigate potential applications and will further strengthen the theoretical analysis presented here.


\section{Appendix}
\label{Sec:Appendix}

\subsection{Proof of Theorem~\ref{Thm:SPGMconv}}
\label{Seq:Proof1}

Consider a single iteration of $\SPGM$ ${\xbm^+ = \Psfhat(\xbm) = \xbm - \gamma\Gsfhat(\xbm)}$, where we used the definition of the gradient mapping. Consider also ${\xbmdetil = \Psf(\xbm) = \xbm - \gamma \Gsf(\xbm)}$. Note that for $B = 1$, we have that
\begin{align*}
&\E[\Psfhat(\xbm)] = \Psf(\xbm) \quad\Rightarrow\quad \E[\Gsfhat(\xbm)] = \Gsf(\xbm)\\
&\E[\|\Psfhat(\xbm)-\Psf(\xbm)\|_2^2]\leq \gamma^2\sigma^2 \quad\Rightarrow\quad \E[\|\Gsfhat(\xbm)-\Gsf(\xbm)\|_2^2]\leq \sigma^2.
\end{align*}
We can then obtain the following bound
\begin{align*}
&f(\xbm^+) = d(\xbm^+) + r(\xbm^+) \\
&\leq f(\xbm) + [\nabla d(\xbm) + \nabla r(\xbmdetil)]^\Tsf(\xbm^+-\xbm) + \frac{L}{2}\|\xbm^+-\xbm\|_2^2 \\
&\quad\quad+ [\nabla r(\xbm^+)-\nabla r(\xbmdetil)]^\Tsf(\xbm^+-\xbmdetil) \\
&\leq f(\xbm) - \gamma \Gsf(\xbm)^\Tsf\Gsfhat(\xbm) + \frac{\gamma^2L}{2}\|\Gsfhat(\xbm)\|_2^2 + \gamma^2 L \|\Gsfhat(\xbm)-\Gsf(\xbm)\|_2^2,
\end{align*}
where the first inequality uses the Lipschitz continuity of $\nabla d$ and twice the convexity of $r$, and the second inequality uses the definition of the gradient mappings, Cauchy-Schwarz inequality, and the Lipschitz continuity of $\nabla r$. By taking the conditional expectation and setting $\gamma = 1/(L\sqrt{T})$, we obtain
\begin{align*}
\E[&f(\xbm^+) - f(\xbm) \,|\, \xbm] \\
&\leq -\gamma \|\Gsf(\xbm)\|_2^2 + \frac{\gamma^2L}{2}(\|\Gsf(\xbm)\|_2^2+\sigma^2) + \gamma^2 L \sigma^2 \\
&\leq -\frac{1}{L\sqrt{T}}\|\Gsf(\xbm)\|_2^2 + \frac{1}{2LT}\|\Gsf(\xbm)\|_2^2 + \frac{3\sigma^2}{2LT} \\
&\leq -\frac{1}{2L\sqrt{T}}\|\Gsf(\xbm)\|_2^2 + \frac{3\sigma^2}{2LT},
\end{align*}
where the final inequality uses the fact that $1/T \leq 1/\sqrt{T}$ for all ${T \geq 1}$. By rearranging the terms and summing up the gradient-mapping norms at different iterations, we finally obtain
\begin{align*}
\E\left[\frac{1}{T}\sum_{t = 1}^T \|\Gsf(\xbm^{t-1})\|_2^2\right] &\leq \frac{1}{\sqrt{T}}\left[2L(f(\xbm^0)-\E[f(\xbm^T)]) + 3\sigma^2\right] \\
&\leq \frac{1}{\sqrt{T}}\left[2L(f(\xbm^0)-f^\ast) + 3\sigma^2\right] ,
\end{align*}
which proves the result.

\subsection{Proof of Theorem~\ref{Thm:signProxconv}}
\label{Seq:Proof2}

Consider a single iteration of $\signProx$ 
$$\xbm^+ = \xbm - \gamma \sgn (\xbm - \Psfhat(\xbm)) = \xbm - \gamma \sgn(\Gsfhat(\xbm)),$$
and a full proximal-gradient iteration ${\xbmdetil = \Psf(\xbm) = \xbm - \gamma \Gsf(\xbm)}$.
We can obtain the following bound
\begin{align*}
&f(\xbm^+) = d(\xbm^+) + r(\xbm^+) \\
&\leq d(\xbm) + \nabla d(\xbm)^\Tsf(\xbm^+-\xbm) + \sum_{i = 1}^n \frac{L_i}{2}(x_i^+-x_i)^2 \\
&\quad\quad + r(\xbm) + \nabla r(\xbm^+)^\Tsf(\xbm^+-\xbm) \\
&= f(\xbm) + [\nabla d(\xbm) + \nabla r(\xbmdetil)]^\Tsf(\xbm^+-\xbm) + \sum_{i = 1}^n \frac{L_i}{2}(x_i^+-x_i)^2 \\
&\quad\quad + [\nabla r(\xbm^+)-\nabla r(\xbmdetil)]^\Tsf (\xbm^+ - \xbm).
\end{align*}
We separately bound the last term above as follows
\begin{align*}
&[\nabla r(\xbm^+)-\nabla r(\xbmdetil)]^\Tsf (\xbm^+ - \xbm) \\
&= [\nabla r(\xbm^+)-\nabla r(\xbm)]^\Tsf (\xbm^+ - \xbm) + [\nabla r(\xbm)-\nabla r(\xbmdetil)]^\Tsf (\xbm^+ - \xbm) \\
&\leq \sum_{i = 1}^n \left[|\nabla r(\xbm^+)_i - \nabla r(\xbm)_i||x_i^+-x_i| + |\nabla r(\xbm)_i - \nabla r(\xbmdetil)_i| |x_i^+-x_i| \right] \\
&\leq \sum_{i = 1}^n [L_i (x_i^+-x_i)^2 + L_i|x_i-\widetilde{x}_i||x_i^+-x_i|] \\
&= \gamma^2 \|\Lbm\|_1 + \gamma \sum_{i = 1}^n L_i|x_i-\widetilde{x}_i| \leq \gamma^2 \|\Lbm\|_1 + \gamma \|\Lbm\|_\infty \|\xbm-\xbmdetil\|_1 \\ 
&\leq \gamma^2 \|\Lbm\|_1 + \gamma^2 \|\Lbm\|_1\|\Gsf(\xbm)\|_1,
\end{align*}
where we used the smoothness assumption on $\nabla r$, the proximal-gradient iterate ${\gamma \Gsf(\xbm) = \xbm - \xbmdetil}$, and the fact that $\|\Lbm\|_\infty \leq \|\Lbm\|_1$. By using this bound in the original inequality, we obtain
\begin{align}
\nonumber&f(\xbm^+) - f(\xbm) \leq -\gamma \Gsf(\xbm)^\Tsf\sgn(\Gsfhat) + \frac{3\gamma^2}{2} \|\Lbm\|_1 + \gamma^2 \|\Lbm\|_1\|\Gsf(\xbm)\|_1 \\
\nonumber&= -\gamma \|\Gsf(\xbm)\|_1 + \frac{3\gamma^2}{2} \|\Lbm\|_1 + \gamma^2 \|\Lbm\|_1\|\Gsf(\xbm)\|_1 \\
\nonumber&\quad\quad + \gamma \Gsf(\xbm)^\Tsf[\sgn(\Gsf(\xbm))-\sgn(\Gsfhat(\xbm))] \\
\label{Eq:ProofBound1}&= -\gamma \|\Gsf(\xbm)\|_1 + \frac{3\gamma^2}{2} \|\Lbm\|_1 + \gamma^2 \|\Lbm\|_1\|\Gsf(\xbm)\|_1 \\
\nonumber&\quad\quad + 2\gamma \sum_{i = 1}^n |\Gsf(\xbm)_i| \indfcn[\sgn(\Gsf(\xbm)_i) \neq \sgn(\Gsfhat(\xbm)_i)],
\end{align}
where $\indfcn[\cdot]$ is an indicator function. The expectation of this function can be further bounded as was done in~\cite{Bernstein.etal2018}
\begin{align*}
\E&[\indfcn[\sgn(\Gsf(\xbm)_i) \neq \sgn(\Gsfhat(\xbm)_i)]] = \P[\sgn(\Gsf(\xbm)_i) \neq \sgn(\Gsfhat(\xbm)_i)] \\
&\leq \P[|\Gsfhat(\xbm)_i-\Gsf(\xbm)_i| \geq |\Gsf_i(\xbm)|] \leq \frac{\E[|\Gsfhat(\xbm)_i-\Gsf(\xbm)_i|]}{|\Gsf(\xbm)_i|} \\
&\leq \frac{\sqrt{\E[(\Gsfhat(\xbm)_i-\Gsf(\xbm)_i)^2]}}{|\Gsf(\xbm)_i|} \leq \frac{\sigma_i}{\sqrt{T}|\Gsf(\xbm)_i|},
\end{align*}
where in the second row we used probability relaxation and the Markov inequality, in the third we used the Jensen's inequality and the variance bound for the mini-batch of size $B = T$. By plugging this expression back into~\eqref{Eq:ProofBound1} and taking the conditional expectation
\begin{align*}
&\E[f(\xbm^+) - f(\xbm) | \xbm] \\
&\leq -\gamma \|\Gsf(\xbm)\|_1 + \frac{3\gamma^2}{2} \|\Lbm\|_1 + \gamma^2 \|\Lbm\|_1\|\Gsf(\xbm)\|_1 + \frac{2\gamma}{\sqrt{T}}\|\sigmabm\|_1 \\
&\leq -\frac{\|\Gsf(\xbm)\|_1}{2\|\Lbm\|_1\sqrt{T}} + \frac{\|\Gsf(\xbm)\|_1}{4\|\Lbm\|_1T} + \frac{3}{8\|\Lbm\|_1T} + \frac{\|\sigmabm\|_1}{\|\Lbm\|_1T} \\
&\leq -\frac{\|\Gsf(\xbm)\|_1}{4\|\Lbm\|_1\sqrt{T}} + \frac{3}{8\|\Lbm\|_1T} + \frac{\|\sigmabm\|_1}{\|\Lbm\|_1T},
\end{align*}
where in the second line we set the step-size to $\gamma = 1/(2\|\Lbm\|_1\sqrt{T})$ and in the last line used the fact that $1/T \leq 1/\sqrt{T}$ and $3/8 \leq 1$. By rearranging the terms and summing up the gradient-mapping norms at different iterations, we finally obtain
$$\E\left[\frac{1}{T}\sum_{t = 1}^T \|\Gsf(\xbm^{t-1})\|_1\right] \leq \frac{4}{\sqrt{T}}\left[\|\Lbm\|_1(f(\xbm^0)-f^\ast) + \|\sigmabm\|_1 + 1\right],$$
which completes the proof.



\newpage

\begin{thebibliography}{10}

\bibitem{Boyd.Vandenberghe2004}
S.~Boyd and L.~Vandenberghe, {\em Convex Optimization}.
\newblock Cambridge Univ. Press, 2004.

\bibitem{Nesterov2004}
Y.~Nesterov, {\em Introductory Lectures on Convex Optimization: A Basic
  Course}.
\newblock Kluwer Academic Publishers, 2004.

\bibitem{Robbins.Monro1951}
H.~Robbins and S.~Monro, ``A stochastic approximation method,'' {\em The Annals
  of Mathematical Statistics}, vol.~22, pp.~400--407, September 1951.

\bibitem{Li.etal2014}
M.~Li, D.~G. Andersen, J.~W. Park, A.~J. Smola, A.~Ahmed, V.~Josifovski,
  J.~Long, E.~J. Shekita, and B.-Y. Su, ``Scaling distributed machine learning
  with the parameter server,'' in {\em Symposium on Operating Systems Design
  and Implementation ({OSDI-14})}, (Broomfield, CO, USA), pp.~583--598, October
  06-08, 2014.

\bibitem{Kamilov.etal2016}
U.~S. Kamilov, I.~N. Papadopoulos, M.~H. Shoreh, A.~Goy, C.~Vonesch, M.~Unser,
  and D.~Psaltis, ``Optical tomographic image reconstruction based on beam
  propagation and sparse regularization,'' {\em IEEE Trans. Comp. Imag.},
  vol.~2, pp.~59--70,, March 2016.

\bibitem{Huang.etal2017}
G.~Huang, Z.~Liu, L.~{van der Maaten}, and K.~Q. Weinberger, ``Densely
  connected convolutional networks,'' in {\em Proc. {IEEE} Conf. Computer
  Vision and Pattern Recognition ({CVPR})}, (Honolulu, HI, USA),
  pp.~2261--2269, July 21-26, 2017.

\bibitem{Seide.etal2014}
F.~Seide, H.~Fu, J.~Droppo, G.~Li, and D.~Yu, ``1-bit stochastic gradient
  descent and its applications to data-parallel distributed training of speech
  {DNN}s,'' in {\em Fifteenth Annual Conference of the International Speech
  Communication Association}, (Singapore), pp.~1058--1062, September 14-18,
  2014.

\bibitem{Alistarh.etal2017}
D.~Alistarh, D.~Grubic, J.~Li, R.~Tomioka, and M.~Vojnovic, ``{QSGD}:
  {C}ommunication-efficient {SGD} via gradient quantization and encoding,'' in
  {\em Proc. Advances in Neural Information Processing Systems 30}, (Long
  Beach, CA, USA), December 4-9, 2017.

\bibitem{Bernstein.etal2018}
J.~Bernstein, Y.-X. Wang, K.~Azizzadenesheli, and A.~Anandkumar, ``{signSGN}:
  {C}ompressed optimization for non-convex problems,'' in {\em Proc. 35th Int.
  Conf. Machine Learning ({ICML})}, (Stockholm, Sweden), July 2018.

\bibitem{Parikh.Boyd2014}
N.~Parikh and S.~Boyd, ``Proximal algorithms,'' {\em Foundations and Trends in
  Optimization}, vol.~1, no.~3, pp.~123--231, 2014.

\bibitem{Moreau1965}
J.~J. Moreau, ``Proximit{\'e} et dualit{\'e} dans un espace hilbertien,'' {\em
  Bull. Soc. Math. France}, vol.~93, pp.~273--299, 1965.

\bibitem{Figueiredo.Nowak2003}
M.~A.~T. Figueiredo and R.~D. Nowak, ``An {EM} algorithm for wavelet-based
  image restoration,'' {\em IEEE Trans. Image Process.}, vol.~12, pp.~906--916,
  August 2003.

\bibitem{Bect.etal2004}
J.~Bect, L.~Blanc-Feraud, G.~Aubert, and A.~Chambolle, ``A $\ell_1$-unified
  variational framework for image restoration,'' in {\em Proc. {ECCV}}
  (Springer, ed.), vol.~3024, (New York), pp.~1--13, 2004.

\bibitem{Daubechies.etal2004}
I.~Daubechies, M.~Defrise, and C.~D. Mol, ``An iterative thresholding algorithm
  for linear inverse problems with a sparsity constraint,'' {\em Commun. Pure
  Appl. Math.}, vol.~57, pp.~1413--1457, November 2004.

\bibitem{Bioucas-Dias.Figueiredo2007}
J.~M. Bioucas-Dias and M.~A.~T. Figueiredo, ``A new {T}w{IST}: {T}wo-step
  iterative shrinkage/thresholding algorithms for image restoration,'' {\em
  IEEE Trans. Image Process.}, vol.~16, pp.~2992--3004, December 2007.

\bibitem{Beck.Teboulle2009}
A.~Beck and M.~Teboulle, ``A fast iterative shrinkage-thresholding algorithm
  for linear inverse problems,'' {\em SIAM J. Imaging Sciences}, vol.~2, no.~1,
  pp.~183--202, 2009.

\bibitem{Eckstein.Bertsekas1992}
J.~Eckstein and D.~P. Bertsekas, ``On the {D}ouglas-{R}achford splitting method
  and the proximal point algorithm for maximal monotone operators,'' {\em
  Mathematical Programming}, vol.~55, pp.~293--318, 1992.

\bibitem{Ng.etal2010}
M.~K. Ng, P.~Weiss, and X.~Yuan, ``Solving constrained total-variation image
  restoration and reconstruction problems via alternating direction methods,''
  {\em SIAM J. Sci. Comput.}, vol.~32, pp.~2710--2736, August 2010.

\bibitem{Boyd.etal2011}
S.~Boyd, N.~Parikh, E.~Chu, B.~Peleato, and J.~Eckstein, ``Distributed
  optimization and statistical learning via the alternating direction method of
  multipliers,'' {\em Foundations and Trends in Machine Learning}, vol.~3,
  no.~1, pp.~1--122, 2011.

\bibitem{Bauschke.etal2008}
H.~H. Bauschke, R.~Goebel, Y.~Lucet, and X.~Wang, ``The proximal average:
  {B}asic theory,'' {\em SIAM J. Optim.}, vol.~19, no.~2, pp.~766--785, 2008.

\bibitem{Yu2013}
Y.~Yu, ``Better approximation and faster algorithm using the proximal
  average,'' in {\em Neural Information Processing Systems ({NIPS})}, (Lake
  Tahoe, CA, USA), pp.~458--466, December 5-10, 2013.

\bibitem{Beck.Teboulle2009b}
A.~Beck and M.~Teboulle, {\em Convex Optimization in Signal Processing and
  Communications}, ch.~Gradient-Based Algorithms with Applications to Signal
  Recovery Problems, pp.~42--88.
\newblock Cambridge, 2009.

\bibitem{Bertsekas2011}
D.~P. Bertsekas, ``Incremental proximal methods for large scale convex
  optimization,'' {\em Math. Program. Ser. B}, vol.~129, pp.~163--195, 2011.

\bibitem{Beck2017}
A.~Beck, {\em First-Order Methods in Optimization}.
\newblock {MOS}-{SIAM} {S}eries on {O}ptimization, SIAM, 2017.

\bibitem{Candes.etal2013}
E.~J. Cand{\`e}s, Y.~C. Eldar, T.~Strohmer, and V.~Voroninski, ``Phase
  retrieval via matrix completion,'' {\em SIAM J. Imaging Sci.}, vol.~6, Feb.
  2013.

\bibitem{Schniter.Rangan2012}
P.~Schniter and S.~Rangan, ``Compressive phase retrieval via generalized
  approximate message passing,'' in {\em Proc. Allerton Conf. on Communication,
  Control, and Computing}, (Monticello, IL), October 2012.

\bibitem{Shechtman.etal2015}
Y.~Shechtman, Y.~C. Eldar, O.~Cohen, H.~N. Chapman, J.~Miao, and M.~Segev,
  ``Phase retrieval with application to optical imaging,'' {\em IEEE Signal
  Process. Mag.}, vol.~32, pp.~87--109, May 2015.

\bibitem{Rudin.etal1992}
L.~I. Rudin, S.~Osher, and E.~Fatemi, ``Nonlinear total variation based noise
  removal algorithms,'' {\em Physica D}, vol.~60, pp.~259--268, November 1992.

\bibitem{Beck.Teboulle2009a}
A.~Beck and M.~Teboulle, ``Fast gradient-based algorithm for constrained total
  variation image denoising and deblurring problems,'' {\em IEEE Trans. Image
  Process.}, vol.~18, pp.~2419--2434, November 2009.

\bibitem{Kamilov2017}
U.~S. Kamilov, ``A parallel proximal algorithm for anisotropic total variation
  minimization,'' {\em IEEE Trans. Image Process.}, vol.~26, pp.~539--548,
  February 2017.

\end{thebibliography}
\end{document}